\documentclass[11pt]{article}

\usepackage{makeidx}
\usepackage{times}      
\usepackage{amsmath}
\usepackage{amsthm}
\usepackage{amssymb}    
\usepackage{graphicx}   
\usepackage{url}
\usepackage[colorlinks]{hyperref}
\usepackage{color}
\usepackage{enumitem}   
\usepackage{booktabs}   
\usepackage{tikz}       
\usepackage{pgfplots}   
\pgfplotsset{compat=1.18}
\usepackage{tikz}
\usepackage{multicol}

\date{}

\begin{document}
\title{A Topology Scavenger Hunt to Introduce Topological Data Analysis}

\author{Lori Ziegelmeier}


\maketitle

\abstract{
Topology at the undergraduate level is often a purely theoretical mathematics course, introducing concepts from point-set topology or possibly algebraic or geometric topology. However, the last two decades have seen an explosion of growth in applied topology and topological data analysis, which are topics that can be presented in an accessible way to undergraduate students and can encourage exciting projects. For the past several years, the Topology course at Macalester College has included content from point-set and algebraic topology, as well as applied topology, culminating in a project chosen by the students. In the course, students work through a topology scavenger hunt as an activity to introduce the ideas and software behind some of the primary tools in topological data analysis, namely, persistent homology and mapper. This scavenger hunt includes a variety of point clouds of varying dimensions, such as an annulus in 2D, a bouquet of loops in 3D, a sphere in 4D, and a torus in 400D. The students’ goal is to analyze each point cloud with a variety of software. This activity can fit nicely into a course where students have been introduced to some of the fundamentals of point-set topology such as connectedness, continuity, compactness, etc. of arbitrary topologies as well as tools from algebraic topology such as simplicial complexes and simplicial homology, which is accessible through the lens of linear algebra. The activity takes approximately a week of class time to provide a brief introduction to persistent homology and mapper, as well as some software resources to perform these computations, and then a week outside of class time for students to work on the scavenger hunt. After completing this activity, students are able to extend the ideas learned in the scavenger hunt to an open-ended capstone project. Examples of past projects include: using persistence to explore the relationship between country development and geography, to analyze congressional voting patterns, and to classify genres of a large corpus of texts by combining with tools from natural language processing and machine learning.}

\section{Introduction}

The field of topology originated in the late 19th century in order to axiomatize and abstract concepts from calculus and analysis \cite{Ziegelmeierstarbird2019topology}. As such, it was considered one of the most theoretical mathematical fields, and courses in topology reflected this structure. However, since the turn of the millenium, there has been an explosion of interest in computational and applied topology to extract the shape of data in modern applications. These applied topology techniques are ripe for explorations which are accessible to students and can be brought into the undergraduate curriculum. 

In a senior-level Topology course, the author introduces some applied topology and topological data analysis (TDA) techniques, such as persistent homology and the mapper graph -- as well as some software tools to compute them -- through a ``Topology Scavenger Hunt". Students explore point clouds of various structures and dimensions in this project, trying to uncover each cloud's topological shape. This activity can fit nicely into a course where students have been introduced to some of the fundamentals of point-set topology as well as concepts from algebraic topology such as simplicial complexes and homology, which are accessible through the lens of linear algebra. While this course has traditionally been taught to a group of 10-20 students, this project could scale to courses of essentially any size and can be accomplished in approximately a week of class time and an additional week for students to work on the scavenger hunt. After completing the Topology Scavenger Hunt, students are well-equipped to explore the structure and shape of data arising from a wide variety of interdisciplinary areas. In fact, students in the course have applied these tools to areas including literature, neuroscience, politics, economics, machine learning, and music, among others.

\section{Prerequisites \& Project Preparation}

Concepts in topology can often be understood without many specific prerequisite courses. However, background in analysis, algebraic structures, and linear algebra can serve as a good foundation and provide necessary mathematical maturity for a student interested in studying topology. While an understanding of concepts from point-set topology is also certainly beneficial, it is not absolutely essential to grasp several applied topology techniques, which find grounding in algebraic and geometric topology. At our institution, the \emph{Topology} course -- which is the only course in topology offered at the college -- first begins by introducing concepts from point-set topology such as topological spaces, continuity, connectivity, compactness, homeomorphisms, and then building new topological spaces through constructions such as disjoint unions, product spaces, or quotient spaces. In recent years, we primarily follow the text \emph{Essential Topology} by Crossley \cite{Ziegelmeiercrossley2006essential}, which covers these topics in a succinct way. Other textbooks including \emph{Introduction to Topology: Pure and Applied} by Adams and Fransoza \cite{Ziegelmeieradams2008introduction}, \emph{Topology} by Munkres \cite{Ziegelmeiermunkres2000topology}, and \emph{Topology Through Inquiry} by Starbird and Su \cite{Ziegelmeierstarbird2019topology} supplement the Crossley text.

\begin{figure}
\centering
\includegraphics[width=.7\textwidth]{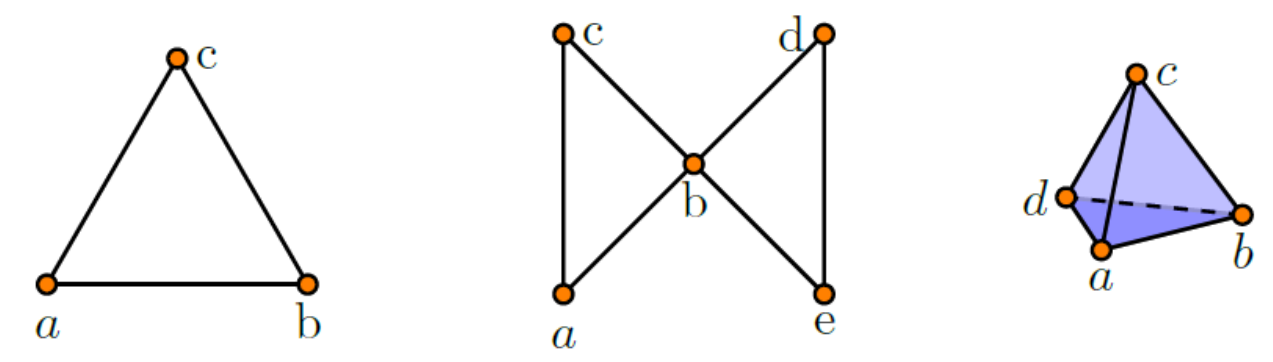} 
\caption{Three simplicial complexes students use to compute persistent homology by hand. (Left) A triangulated circle $S^1$, which has Betti numbers $\beta_0=1,\beta_1=1,\beta_2=0$; (Middle) a triangulated figure 8, which has $\beta_0=1,\beta_1=2,\beta_2=0$; (Right) a triangulated sphere, which has $\beta_0=1,\beta_1=0,\beta_2=1$.}
\label{fig:Ziegelmeiertriangulatedsimpcomp}
\end{figure}

Once the foundation of these topology underpinnings has been laid, we then move on to combinatorial and algebraic topology concepts such as simplicial complexes and simplicial homology which can serve as the primary building blocks of applied topology techniques. For the reader unfamiliar with these concepts, the Hatcher \emph{Algebraic Topology} text \cite{ZiegelmeierHatcher} is a classic resource. These concepts are also discussed in the Crossley text \cite{Ziegelmeiercrossley2006essential}, and a fun resource is ``Self-help homology tutorial for the simple(x)-minded" by Topaz \cite{Ziegelmeierchadtopaz_2021}. We briefly introduce these concepts in Appendix \ref{app:ZiegelmeierDefinitions}. Once these definitions have been established, students work through explicit examples of computing bases and dimensions for the boundaries, cycles, and homology classes of various simplicial complexes, such as those in Figure \ref{fig:Ziegelmeiertriangulatedsimpcomp}, by hand.
This can be done with standard linear algebra techniques to compute images and kernels of boundary matrices. After these activities, students perform homology computations using the software Macaulay2 \cite{ZiegelmeierMacaulay2} to explore larger simplicial complexes such as a triangulated torus or Klein bottle.

After defining key terms and working carefully through examples in combinatorial and algebraic topology, students are now equipped to understand concepts in applied topology. To provide some background material, students read the article ``Barcodes: the persistent topology of data" by Ghrist \cite{ZiegelmeierbarcodeGhrist} and are assigned guiding questions to scaffold their reading (included in Appendix \ref{ZiegelmeierReading} with example student solutions in Appendix \ref{ZiegelmeierReadingSolutions}). Optionally, students are also encouraged to read articles such as \cite{Ziegelmeiercipra2009,ZiegelmeierMunch_2017,ZiegelmeierLesnick}, or the ``Topological Data Analysis and Persistent Homology" section of \cite{ZiegelmeierTZH15}. Students are also directed to the Applied Algebraic Topology Research Network (AATRN) YouTube channel \cite{ZiegelmeierAATRN}, which has several excellent applied topology research talks and tutorials. More comprehensive references for the instructor include "Topology and Data" by Carlsson \cite{ZiegelmeierCarlsson2009TopologyAD}, "Persistent homology--A Survey" by Edelsbrunner and Harer,  \cite{Ziegelmeierpersistenthomologyasurvey}, \emph{Elementary Applied Topology} by Ghrist \cite{Ziegelmeierghrist2014elementary}, \emph{Computational Topology - an Introduction} by Edelsbrunner and Harer \cite{ZiegelmeierEdeHar2010}, and \emph{Computational Topology for Data Analysis} by Dey and Wang \cite{Ziegelmeierdey2022computational}.

The two most common applied topology techniques to extract the shape and structure from data are \emph{persistent homology}, which measures homology of a sequence of nested simplicial complexes defined using a scale parameter to track which features persist through these scales, and the \emph{mapper graph}, a 1-dimensional structure that is useful for exploration and visualization of the data. We briefly introduce these concepts in Appendix \ref{app:ZiegelmeierTDADefinitions}.

Persistent homology can be computed using linear algebra techniques. However, when working with real data, simplicial complexes can result in boundary matrices with millions or billions of columns. As such, there are a number of research groups working in computational topology to investigate and develop efficient algorithms to perform these computations on large-scale data. Thus, there are a variety of software packages available to compute both persistent homology and the mapper graph, with more being developed regularly. In our course, we direct students to the R TDA package \cite{Ziegelmeierfasy2015introduction}, Ripser \cite{ZiegelmeierBauer2021Ripser} (including the online browser version at \cite{ZiegelmeierBrowserRipser}), Eirene \cite{Ziegelmeierhenselmanghristl6}, Python scikit-tda \cite{Ziegelmeierscikittda2019}, and giotto-tda \cite{Ziegelmeiergiotto} for persistent homology and the TDAmapper package in R \cite{ZiegelmeierTDAmapper} to compute the mapper graph. Several of these packages have associated tutorials and demonstrations, and the author of this article, with co-organizers Matthew Wright and Matthew Richey, ran a CBMS regional research conference in 2017 which centralized a few of these resources in the following website \cite{Ziegelmeierwright5_2017}. In fact, the point clouds used in the Topology Scavenger Hunt were developed by Wright and the author for that conference.

\section{Course Project}

After spending approximately a week of class time\footnote{This does not include the time devoted to discussing the more standard algebraic topology methods of simplicial complexes and simplicial homology, which typically takes an additional week or two of class time.} introducing students to persistent homology (primary focus) and the mapper graph (cursory discussion) as well as a brief overview of some of the software packages available to perform these computations, students are now equipped to begin the Topology Scavenger Hunt project. The directions for the project are succinct and as follows:

\begin{quote}
    Your task is to complete a topology scavenger hunt by analyzing the point cloud data posted in [course management site].
There are point clouds from 2D, 3D, 4D, and 400D. Your task is to take each one and determine the
topological structure of the point cloud. Use at least three different pieces of software to do so (one
software per point cloud is fine, but make sure to try a variety of approaches on different software for
different point clouds). Of course, you can plot the 2D and 3D data sets to verify your intuition (which
I encourage you to do \textbf{after} you have first tried to determine the structure). Explain your approach,
submit your code with any relevant plots, and include your guess for the structure of each point cloud.
\end{quote}

Most students in recent iterations of the course gravitate toward using the R TDA package and the R TDAmapper package -- likely due to familiarity with using R in our introductory statistics and linear algebra courses -- and the online browser version of the software Ripser, which is both fast and straightforward. Students are given approximately a week of time outside of class to complete the project.

\begin{figure}
\centering
\includegraphics[width=\textwidth]{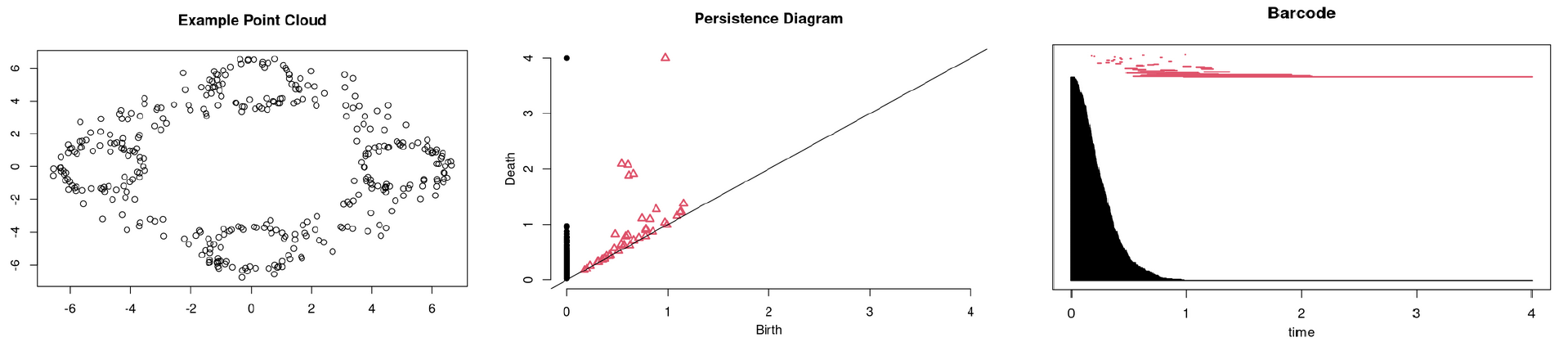} 
\caption{Example of one of the point clouds (Left) to be analyzed in the Topology Scavenger Hunt as well as the associated persistence diagram (Middle) and barcode (Right) computed using the R TDA package \cite{Ziegelmeierfasy2015introduction}. The persistence diagram and barcode display the same information in different ways. The 0 and 1-dimensional features are computed and displayed in black and red, respectively.}
\label{fig:ZiegelmeierpointsPDBarcode}
\end{figure}

The point clouds for the scavenger hunt can be found at the following GitHub page \cite{Ziegelmeierpointclouds}. We display one such point cloud (points2.csv) in Figure \ref{fig:ZiegelmeierpointsPDBarcode} as well as the 0- and 1-dimensional persistence diagrams and barcodes. Ideally, a student would begin exploration without actually plotting the point cloud first and would instead use the persistence diagram and barcode to topologically infer the structure of the point cloud as now described. 

Recall that 0-dimensional features (represented in black in the persistence diagram and barcode) relate to the number of connected components, and 1-dimensional features (represented in red) relate to the number of loops. Since there is one black isolated point in the persistence diagram and one long bar in the barcode, we infer that all of the points join into a single connected component fairly early (i.e., at $\varepsilon\approx 1$) in the filtration. In the 1-dimensional signatures, there is one red triangle that is far away from the diagonal and one long red bar, indicating a 1-dimensional topological loop that is born when all points are connected into a single component and persists through the entirety of the displayed filtration. There are also four other red triangles further from the diagonal and corresponding four bars of longer length than the majority of features but of shorter persistence than the red isolated point in the persistence diagram or long bar in the barcode. These can be interpreted as four topological loops that are born before all of the points connect into a single component and persist until $\varepsilon\approx 2$. That is, they are significant but shorter-lived 1-dimensional features. All other red features are either close to the diagonal in the persistence diagram or short bars in the barcode and can be intuited as noise. Putting these pieces together, we can infer that the point cloud is primarily a single structure, with one larger topological loop and four smaller loops. Note that homology cannot describe precisely where these loops are located in relation to one another, however exploration with techniques such as the mapper graph can further illuminate the loops' locations. 

Of course, the power of these TDA techniques is to explore high dimensional data which cannot simply be visualized by plotting a point cloud as in the case of Figure \ref{fig:ZiegelmeierpointsPDBarcode}. However, exploration with these lower-dimensional, pedagogical point clouds can help students refine and expand their understanding of these tools. The goal of the Topology Scavenger Hunt is to equip students with the topological tools necessary to explore the shape and structure of a variety of data sets. Then, they are able to take these tools to go explore their own projects on interdisciplinary topics which appeal to them. The creativity of these interdisciplinary projects comes from the students.

\section{Extensions \& Alterations} \label{Ziegelmeiersec:extensions}

Students in the \emph{Topology} course take the tools learned throughout the class and then go on to explore an open-ended capstone project selected by the student in the last few weeks of the semester. Sometimes this involves theory related to topology concepts, for example, exploring more in-depth separability theorems, cohomology, the classification theorem of compact surfaces, or Sperner's Lemma. More often than not, though, students are inspired by the Topology Scavenger Hunt to explore concepts in applied topology to analyze a data set of interest arising from an interdisciplinary application, to understand in detail the algorithms used to compute persistence, or to learn more about a concept not discussed in class (a few topics are mentioned in the last paragraph of this section).

There have been a number of interesting capstone projects related to interdisciplinary applications of these tools. A few examples are now highlighted: 

\begin{itemize}
    \item One student used TDA combined with natural language processing to analyze hundreds of books in order to differentiate them by their genre. He used persistent homology to extract the shape of each text and then performed classification of these topological signatures using a t-test on homological statistics as well as a clustering algorithm to confirm that TDA provides a powerful way of analyzing text. 
    \item Another student -- expanding upon the study of \cite{ZiegelmeierPiangerelli} -- used a topological measure called persistent entropy to differentiate between EEGs of patients with seizures versus no epileptic activity as well as EEGs which have been labeled ``normal" versus ``abnormal." In contrast to the previous work which analyzed hour-long EEGs, the student found that persistent entropy was not overwhelmingly effective in classifying short (20-second) segments of EEGs, although segments of focal seizures did tend to have higher entropy than background activity for each individual patient. \item A group of students analyzed persistent homology cycle representatives arising from analysis of Hegel's ``Phenomenology of Spirit" in order to explore a correspondence between loops containing a set of significant words and their meaning. 
    \item Two other students combined their interests in geography with topology to analyze cities across the United States to investigate questions concerning racial segregation and food deserts.
    \item Another group applied TDA to fragments of European classical music by mapping musical events in scores to sets of points in Euclidean space, performing persistent homology, then using classification techniques to reveal stylistic differences and comparisons of composers.
    \item Another student used persistence to explore the relationship between country development and geography by considering various indicators such as gross domestic product per capita, average life expectancy, infant mortality, etc. to quantify the development of each country based on persistent homology clusters and cycles. This project was developed beyond the capstone course and turned into a publication \cite{ZiegelmeierBanZie2018}. 
\end{itemize}

\paragraph{TDA and Polarized Congressional Voting} To expound upon the possibilities of analyzing an interdisciplinary application using TDA, we now describe an additional project completed by a student in more detail. This student used TDA to analyze polarization in United States congressional voting patterns over time. In particular, he computed the ideological distance between congressional members of each congress (split into the House and Senate) based on DW-NOMINATE scores \cite{ZiegelmeierPoole, ZiegelmeierNOMINATE}, a method for turning high-dimensional roll call votes into a two-dimensional ideological score. He then computed a pairwise Euclidean distance matrix to capture the (dis)similarity among all congressional members for each congress (up to the 116th congress from 2019-2020). Given such a distance matrix, persistent homology can be used to measure the shape of the corresponding data. As such, he computed distance matrices and persistent homology for the House and Senate of each congress. After visually looking through barcodes of various congresses, he noticed differences between the shape of polarized versus unpolarized congresses. As such, he computed the bottleneck distance \cite{ZiegelmeierEdeHar2010} between each persistence barcode corresponding to a congress and compared this to the absolute difference in polarization (defined as the absolute difference of the average democrat and the average republican DW-NOMINATE scores) via a linear regression model.\footnote{The student computed the regression from 1880 onward, as this provided two stable political parties throughout the sample.} The relationship between the bottleneck distance and the absolute difference in polarization was both positive and significant, showing that congress sessions with similar polarization scores are also similar in these topological measures.  

\begin{figure}
\centering
\includegraphics[width=\textwidth]{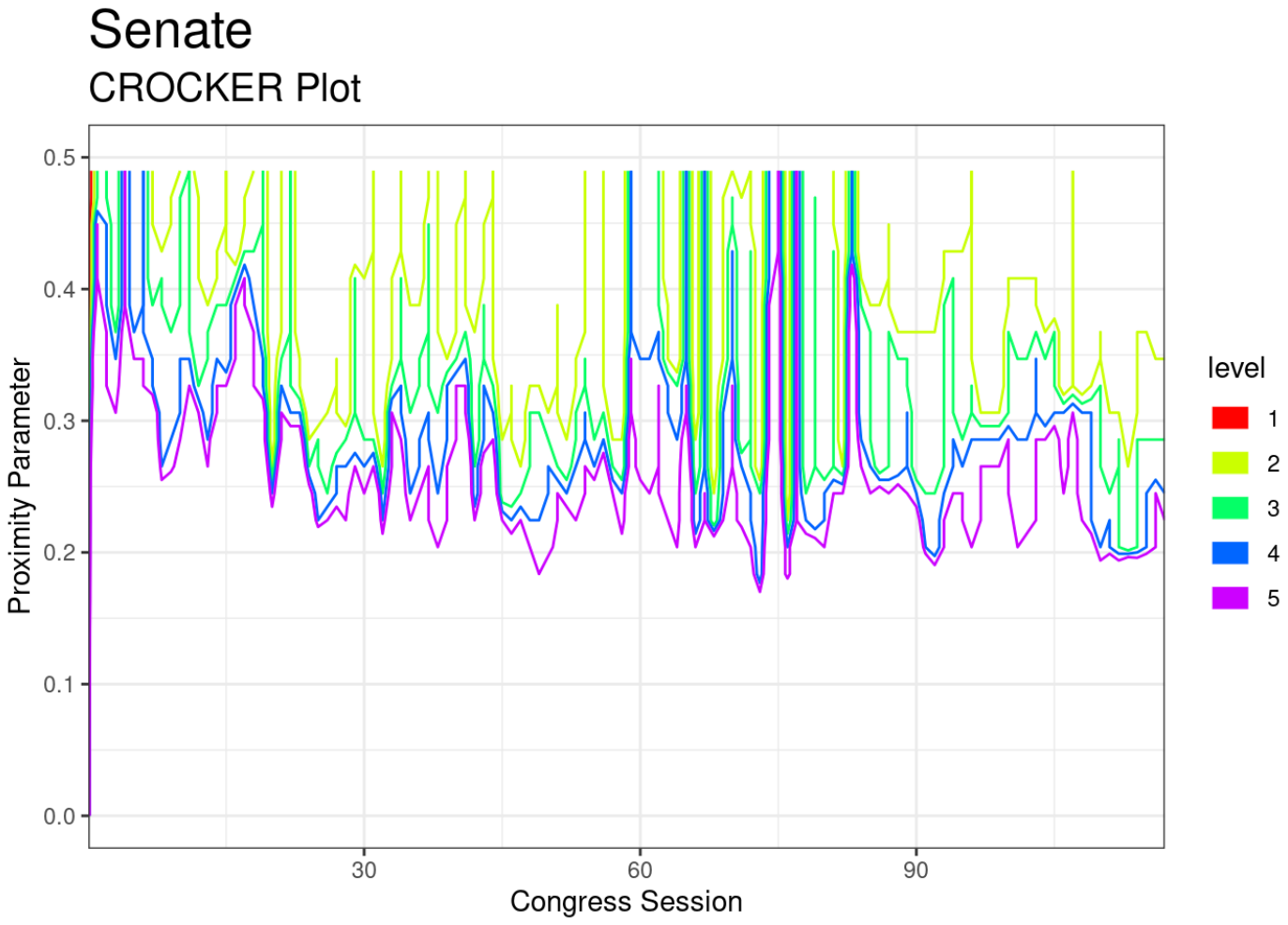} 
\caption{Contour plot of $\beta_0(t,\varepsilon)$, where $t$ corresponds to congressional session and $\varepsilon$ is topological proximity or scale parameter, for the United States Senate. Five contour levels are displayed, interpreting $\beta_0 > 5$ as noise.}
\label{fig:ZiegelmeierSenateCrocker}
\end{figure}
Further, this student constructed a CROCKER plot, a visualization and vectorization of topological information for multiple input parameters described in the article \cite{ZiegelmeierTZH15}, for the House and for the Senate through all congressional sessions up to the 116th. For this project, the CROCKER can be thought of as a function giving the Betti number $\beta_k(t,\varepsilon)$ of two inputs, congressional session $t$ and topological proximity or scale parameter $\varepsilon$. It can be constructed as a matrix and visualized as a contour diagram; see Figure \ref{fig:ZiegelmeierSenateCrocker} for an example, which displays the $\beta_0$ topological information of the Senate. Recall that $\beta_0$ refers to the number of connected components within a simplicial complex. Thus, the value of $\beta_0(t,\varepsilon)$ reveals the number of connected components of a particular congress at a given topological scale. In the region of the 73rd-84th congresses (1934-1956), we observe that three, four, and five connected components persist for a much broader range of scale values than the surrounding congresses. This is a region corresponding to low polarization according to the DW-NOMINATE score. In contrast, regions of high polarization---e.g., the 45th-58th congresses (from 1879-1905) and more recent years---have markedly different structure revealed by the CROCKER. For instance, for the region corresponding to the 89th-92nd congresses, we see that there is greater separation between the second and third contour lines, but all senators connect into a single component at a smaller proximity parameter. 

A plausible interpretation of this follows: In unpolarized eras, there was more variety (reflected in greater proximity values and more connected components) of ideology while in more polarized eras, two connected components at smaller scales are a more central feature, reflecting the more dramatic polarization between the two-party system of the United States Senate. One hypothesis for this behavior is that even though political polarization is defined as a large difference between the party means, there could be less ideological variation overall, but in unpolarized eras, there could be a few members whose voting patterns are not similar to anyone else's. This example highlights how TDA can provide more nuanced information than the more traditional method of simply examining means. More analysis into which senators make up each connected component could be done to better understand these connections. A data set such as this is ripe for more questions and exploration than can be feasibly completed in a couple week course project, and further analysis could lead to interesting insights.

\paragraph{Extensions and Alterations} This student project provides an example of interdisciplinary explorations students can achieve after completing the Topology Scavenger Hunt project. The goal of this activity is to introduce students to tools for them to explore something of interest in order to take ownership of their capstone project. In addition, some students have been inspired to continue learning about TDA beyond the course by doing summer research and independent studies, which have led to fruitful collaborations. While not every collaboration has resulted in publishable work, some have resulted in the following articles \cite{ZiegelmeierUlmZieTop2019, ZiegelmeierXiaAdaTop2021, ZiegelmeierLiThoHen2021optimal}, corresponding to analyzing models of pea aphid movement using more traditional order parameters and topology, an extension of the CROCKER tool discussed above as well as a theoretical exploration and clustering of simulations from a well-known biological aggregation model, and an empirical investigation of a variety of data sets on different notions of minimal cycle representatives, respectively.

There are a number of extensions and alterations that can be made to the scavenger hunt and then the subsequent capstone project. Of course, a wide variety of point clouds could be used, including real-world data. Alternatively, there are notions of persistent homology in a variety of contexts including: (sub/super) level-set filtrations (to analyze functions or surfaces) \cite{ZiegelmeierEdeHar2010}, cubical complexes (to analyze image or voxel data) \cite{Ziegelmeierbleile2021persistent}, zig-zag modules (to consider not just a sequence of increasing filtrations but a 1-parameter family with maps in either direction) \cite{ZiegelmeierCardeSMor2009}, sliding windows of a signal (to analyze time series) \cite{Ziegelmeierwindowsandpersistence}, or multidimensional (more recently referred to as multiparameter) persistence \cite{Ziegelmeiermultidimensional}. Further, there are a number of approaches to compare persistence barcodes and diagrams including: the bottleneck or Wasserstein distances \cite{ZiegelmeierEdeHar2010}, means of a collection of diagrams \cite{ZiegelmeierTurMilMuk2014}, confidence intervals \cite{Ziegelmeierfasyconfidence}, or encoding as a vectorization to be used in machine learning \cite{Ziegelmeierbubenik2015statistical, Ziegelmeierreininghaus2015stable, ZiegelmeierAdaEmeKir2017,Ziegelmeierjournals/cgf/CarriereOO15}. There is also a vast amount of theory related to applied topology with perhaps the most central being the stability of persistence, e.g., that the distance between persistence diagrams of two persistence modules is bounded by the distance between the two input functions \cite{ZiegelmeierstabilityCohen} or point clouds \cite{Ziegelmeierccgmo-ghsssp-09}.

\section{Reflections \& Advice}

The key to having the Topology Scavenger Hunt project succeed is to carefully motivate underlying concepts in class; that is, thorough discussions of simplicial complexes, homology, and persistent homology. The author does this through in-class annotated notes, showing slide presentations of concepts and how they relate to applications, and in-class student activities. Once this foundation is laid, the other primary challenge is to make sure that students have access to a variety of working software. Mathematics students at our college generally have a good grasp on computing, being introduced to a variety of software in our introductory computer science (Python), statistics (R), calculus (R and Mathematica), and linear algebra (R) courses as well as many of our mid- and upper-level mathematics courses (including R, MATLAB, Python, Mathematica, SAGE, Octave). If students do not have this previous exposure, more time may need to be spent on introducing the software. Of course, it is infeasible for the instructor to know how to implement all of the TDA software available but having a couple of demonstrations developed is key. As is the case with any software, it can sometimes be challenging to make sure it is installed properly with appropriate dependencies, and so verifying that the software is working as expected, well in advance of presenting it, is essential.

The main goal of the author for this project is that students explore TDA techniques and software and have fun while doing so. Generally, grading of the Topology Scavenger Hunt project is quite ``loose" in that students who have attempted to analyze most of the point clouds with at least three different types of software and generally make correct inferences about the structure of each point cloud receive full credit. Then, students can take these tools to explore interdisciplinary applications which appeal to them. Through the Topology Scavenger Hunt, students get to see that what once was a very theoretical field in topology now, with the power of TDA, can be used to explore the shape of a variety of interdisciplinary data.

\section*{Acknowledgement}

The author would like to acknowledge and thank two former students, Aidan Toner-Rodgers and Ari Holcombe-Pomerance. The former student's topology capstone project was discussed in Section \ref{Ziegelmeiersec:extensions} and the latter's solutions appear in Appendix \ref{ZiegelmeierReadingSolutions}. Ziegelmeier was supported by NSF grant CDS\&E-MSS-1854703 during development of this course project.

\bibliographystyle{abbrv}
\bibliography{lz_bib}

\appendix 
\section*{Appendix}

\section{Algebraic Topology Definitions}\label{app:ZiegelmeierDefinitions}

In this appendix, we provide an overview of algebraic topology definitions fundamental to tools in topological data analysis.

\paragraph{Simplicial Complexes} A \emph{simplicial complex} is a way to turn a topological space into a combinatorial object, a higher-dimensional analog of a network or graph structure, by decomposing it into simple pieces called \emph{simplices}. Each of these $k$-simplices are the smallest convex set containing $k+1$ points e.g., a 0-simplex is a vertex, a 1-simplex is an edge, a 2-simplex is a filled-in triangle, a 3-simplex is a solid tetrahedron, and so on. The collection of simplices satisfies two properties (1) two $k$-simplices are either disjoint or intersect in a lower-dimensional simplex, and (2) any subset in a simplicial complex must contain all lower-dimensional simplices. Figure \ref{fig:Ziegelmeiertriangulatedsimpcomp} displays examples of three different simplicial complexes.

\paragraph{Boundaries and Cycles} For each dimension $k\geq 0$, we create an abstract vector space $C_k$ with basis consisting of the set of $k$-simplices within a simplicial complex and coefficients (typically) coming from a field. Elements of $C_k$ are called \emph{$k$-chains}. We algebraically describe the boundary of an ordered $k$-simplex $(v_0,v_1,\ldots,v_k)$ by the following \emph{boundary map} $\partial_k:C_k \to C_{k-1}$,
$$
        \partial_k(v_0, v_1, \ldots, v_k) = \sum_{i=0}^k (-1)^i (v_0, \ldots, \hat{v_i}, \ldots, v_k),
$$
where $\hat{v_i}$ omits node $v_i$. The collection of these vector spaces $C_k$ along with the associated boundary maps $\partial_k$ form a \emph{chain complex}
\[\ldots C_{k+1} \xrightarrow{\partial_{k+1}} C_{k} \xrightarrow{\partial_{k}} C_{k-1} \xrightarrow{\partial_{k-1}} \ldots \xrightarrow{\partial_3} C_2 \xrightarrow{\partial_2} C_1 \xrightarrow{\partial_1} C_0 \xrightarrow{\partial_0} 0. \]
Within each vector space $C_k$, there are two important subspaces (1) the set of \emph{boundaries}  given by the image of the boundary map $B_k:=\textbf{im}( \partial_{k+1})$, and (2) the set of \emph{cycles}, which vanish under the boundary map, e.g., elements in the kernel $Z_k:=\textbf{ker}( \partial_k)$. The boundary map satisfies the following fundamental property $\partial_k \circ \partial_{k+1}=0$, so the set of boundaries is always contained inside the set of cycles.

\paragraph{Homology} \emph{Homology} uncovers those cycles which are not also boundaries by computing the quotient $H_k:=Z_k/B_k$. This splits the cycles into equivalence classes, where two cycles are \emph{homologous} (equivalent) if their difference is a boundary of a ($k+1$)-chain. An element of $H_k$ is called a \emph{$k$-dimensional homology class}, which can be given by a particular \emph{cycle representative}. The dimension of the homology vector space is called the \emph{Betti number}, $\beta_k:=\textbf{dim}(H_k)=\textbf{dim}(Z_k)-\textbf{dim}(B_k)$, which counts the number of independent holes of dimension $k$. That is, a hole surrounded by $k$-simplices. For instance, $\beta_0$ counts the number of connected components in the simplicial complex, $\beta_1$ counts the number of topological loops that bound a 2D void, $\beta_2$ counts the number of trapped 3D volumes, and so on as dimension increases. Betti numbers are topological invariants, meaning that they are unchanged under continuous deformations of the object such as stretching, compressing, warping, and bending. Thus, they measure inherent qualitative properties about the shape of an object.

\section{Reading Questions} \label{ZiegelmeierReading}

Your main assignment for this week is to read the article ``Barcodes: the Persistent Topology of Data" by
Robert Ghrist \cite{ZiegelmeierbarcodeGhrist}. To guide you in your reading, please answer the following questions:
\begin{enumerate}
	\item What is a Cech complex? What is a Rips complex? Construct a toy example that illustrates when these two complexes have a different structure from one another.
	\item What is the optimal $\varepsilon$ that should be chosen to construct the Rips or Cech complex? Hint: This is a
trick question! Why?
\item Section 2.2 gets a bit technical. What is the main point of this section?
\item What is a Betti number? Why is it useful? How is it displayed in a barcode?
\item Section 3 describes a really cool application to natural images. Describe the process that Mumford and
others undertook to analyze natural images. Here are some questions to guide you, but fill in some of
the details.
\begin{enumerate}
	\item Why is the original data set on a topological seven-sphere $S^7$?
	\item What is going on with the codensity function?
	\item Why does a different choice of density and neighbor parameter produce different barcodes representing homology?
	\item Why are there 5 $H_1$ generators when $T = 25$ and $k = 15$?
	\item How does a Klein bottle seem to appear in this example?
\end{enumerate}
\item Provide a 1-2 paragraph summary of the article.
\end{enumerate}

\section{Reading Questions Possible Solutions} \label{ZiegelmeierReadingSolutions}

Answers to the previous questions will likely vary quite a bit across students. Here are some solutions from a student in a recent iteration of the course, with minimal editing by the professor:

\begin{enumerate}
	\item A Cech complex is a simplicial complex where the $k$-simplices are determined by ($k+1$)-tuples of points in a Euclidean space whose neighborhoods intersect.
		(For example, vertices are determined by points, edges by 2 points with intersecting neighborhoods, triangles by 3 points with intersecting neighborhoods, etc.).
		A Rips complex is similar, but instead of determining tuples by intersecting neighborhoods, we let $k$-simplices correspond to ($k+1$)-tuples within a certain pairwise distance.
		These are generally different structures, as the Rips complex generally has more simplices (since we only consider distance rather than entire neighborhoods).
	\item There is no optimal choice for $\varepsilon$!
		We might want to increase $\varepsilon$ to remove small holes between relatively close data points, but in doing so we risk removing more significant holes in our data set.
		The goal is then to strike a balance between too much and too little detail. Therefore, persistent homology tracks topology through \textbf{all} choices of $\varepsilon$ to see which features persist across the filtration.
	\item This section gives a definition for persistent homology (it is the image of an induced homomorphism).
		The main point is to detail some of the algebra underlying the creation of barcodes, which are introduced in the following section.
	\item The Betti number is the rank of $H_k$ (the $k$-th homology group).
		These are useful because they tell us about the structure of our simplicial complex.
		(For example, $\beta_0$ is the number of connected components, $\beta_1$ is the number of holes, $\beta_2$ is the number of trapped volumes, etc.).
		Ghrist calls barcodes the ``persistence analogue'' of Betti numbers, since they convey information of rank through the filtration by looking at vertical slices at a given value of $\varepsilon$ and counting the number of bars that intersect this slice.
	\item Mumford et al. start with 5000 3x3 pixel squares chosen randomly from an image and store these in $\mathbb{E}^9$.
		By normalizing and selecting only high-contrast images, they reduce the dimensionality and project the data onto the 7-sphere.
		The codensity function is used to evaluate the closeness of neighbors to a given point (the density of the dataset in that area).
		It is used to filter data by creating ``cores'' for point clouds.
		Now that we have a way to homologize the data set, we can examine its persistent homology (for example, $H_1$ tells us about ``loops'' in the data set).
		We can change the codensity parameter $k$ to get different point clouds resulting in different persistent homology for each.
		A lower value of $k$ generates barcodes that place more emphasis on local densities rather than averaging overall density.
        There are 5 $H_1$ generators resulting from a primary circle and two secondary circles which intersect at common 3x3 patches.
		Given a specific choice of $k$ and $T$, the authors were able to achieve the Klein bottle from this process, illustrated with an identification polygon with two opposite edges having the same orientation and the other two having opposite orientations.
	\item In ``Barcodes: The persistent topology of data,'' Robert Ghrist introduces an application of persistent homology onto point-cloud data sets by creating a topological characterization called a barcode.
		Barcodes work in a similar way to Betti numbers, by analyzing simplicial complexes to give an overview of a shape's structure across varying scale parameter $\varepsilon$.
		We can find the barcode of high-dimensional data as a means of shape recognition to gain insight about a data set's structure.
		One such application, introduced by Mumford et al. is to characterize the shape of pixel patches in natural images.
\end{enumerate}

\section{Topological Data Analysis Definitions} \label{app:ZiegelmeierTDADefinitions}

This appendix describes key definitions for two common topological data analysis techniques: persistent homology and the mapper graph.

\paragraph{Filtrations} A common method to form a simplicial complex from a set of discrete data equipped with a notion of distance is called a \emph{Vietoris-Rips complex} (alternatively, \emph{flag} or \emph{clique} complex). One first defines a \emph{scale parameter} (alternatively, \emph{proximity} or \emph{filtration} parameter) $\varepsilon$ at which connections between proximate nodes are made. For every collection of $k+1$ nodes which are pairwise less than $\varepsilon$, a $k$-simplex is formed. For instance, if four nodes are all pairwise within scale $\varepsilon$, then we connect all six edges, fill in each of the four triangles bounded by those edges, and fill in the solid tetrahedron bounded by the four triangles to get a $3$-simplex. 

\paragraph{Persistent Homology} In practice, it is often challenging to determine a single appropriate scale at which these connections are made, particularly without \emph{a priori} knowlege. Instead, a nested sequence of simplicial complexes, called a \emph{filtration}, $S_{\varepsilon_0} \subseteq S_{\varepsilon_1} \subseteq \cdots \subseteq S_{\varepsilon_j} \subseteq \cdots$ indexed by increasing scale parameters $\varepsilon_0 \leq \varepsilon_1 \leq \cdots \leq \varepsilon_j \leq \cdots$ is formed. To each simplicial complex, we apply $k$-dimensional homology $H_k$. The inclusion of simplicial complexes through the filtration induces linear maps on homology $$H_k(S_{\varepsilon_0}) \to H_k(S_{\varepsilon_1}) \to \cdots \to H_k(S_{\varepsilon_j}) \to \cdots .$$ 
This sequence of spaces and associated maps is known as a  \emph{persistence module}. \emph{Persistent homology} then tracks classes of each homology space through the filtration, often referred to as topological \emph{features}. The word persistence corresponds to the ranges of the scale parameter over which features persist, with short-lived features often intuited as noise and long-lived as true signal. The scale at which a feature (hole) first appears is given by its \emph{birth} value and the scale at which a feature disappears is its \emph{death} value. A feature remaining throughout the entirety of the filtration has death at infinity. The \emph{persistence} (or \emph{lifetime}) of a feature is then its death minus birth. 

\begin{figure}
\centering
\includegraphics[width=\textwidth]{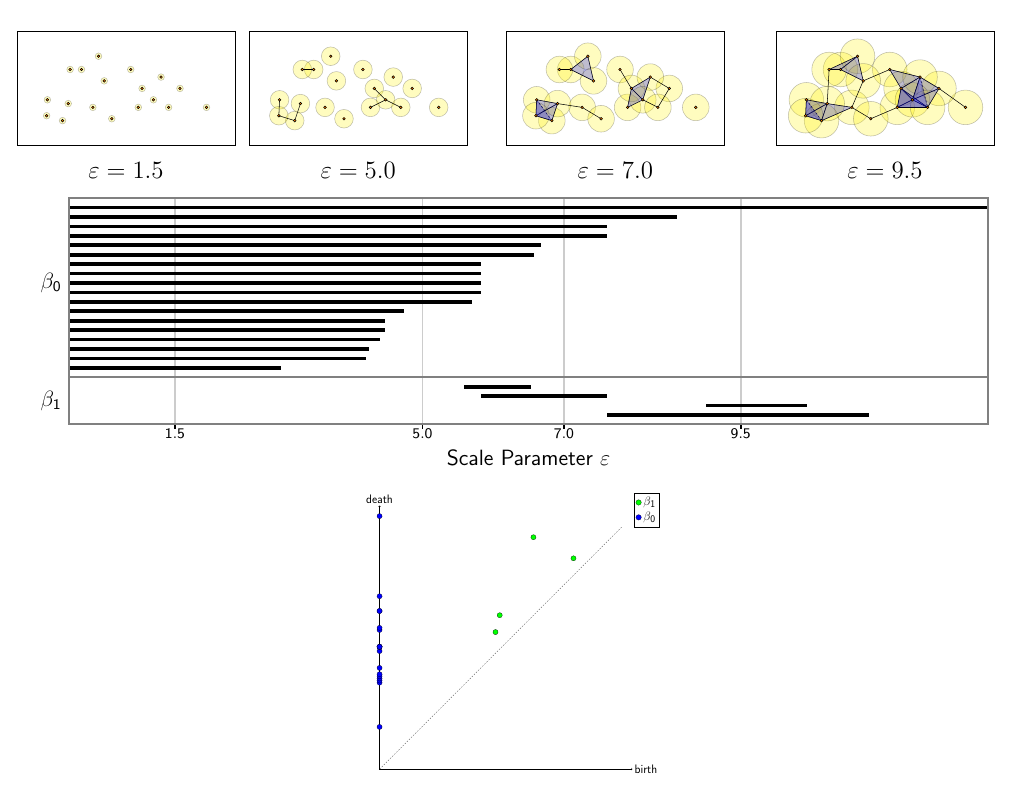} 
\caption{Example from the author's article \cite{ZiegelmeierTZH15}, available under Creative Commons License \url{https://creativecommons.org/licenses/by/4.0/}, of the topological barcode (Middle) and persistence diagram (Bottom) of a Vietoris-Rips complex. The top four figures display the simplicial complex of 18 points for different values
of the scale parameter $\varepsilon$, where the yellow circles have radius $\varepsilon/2$ and are used for visualization but are not actually part of the simplicial complex. The vertical lines in the barcode correspond to these four levels of $\varepsilon$. The number of horizontal bars intersecting each line gives
the values of $\beta_0=18$ and $\beta_1=0$ at $\varepsilon=1.5$, $\beta_0=11$ and $\beta_1=0$ at $\varepsilon=5.0$, $\beta_0=4$ and $\beta_1=1$ at $\varepsilon=7.0$, and $\beta_0=1$ and $\beta_1=2$ at $\varepsilon=9.5$. Each bar in the barcode has an equivalent persistence point in the persistence diagram with coordinates given by its birth and death values. In a Vietoris-Rips construction, the 0-dimensional features will all be born at $\varepsilon=0$ and join into a single component once the scale parameter is large enough.}
\label{fig:Ziegelmeierpedexample}
\end{figure}

\paragraph{Visualization of Topological Features} Enumeration of topological features through a filtration as well as their associated birth and death scales is often encoded in graphical representations known as \emph{barcodes} or \emph{persistence diagrams}. In the former, there is a distinct interval (or bar) corresponding to the range of scales over which each homology class persists. That is, if a feature is born at scale $b$ and no longer remains at scale $d$, then there is an associated half-open interval $[b,d)$ representing the lifetime of that feature. The number of $k$-dimensional holes given by the Betti number $\beta_k$ at scale $\varepsilon$ is the number of distinct bars intersecting the vertical line through $\varepsilon$. An equivalent topological representation is a persistence diagram, which encodes each feature as a multi-set of points in the plane, where the $x$-axis indicates the birth coordinate and the $y$-axis is the death. Thus, a feature (hole) is represented as a point $(b,d)$ in a persistence diagram. All features must be born before they die, so these points lie on or above the diagonal. Points that are farther from the diagonal have long persistence, while short-lived features are near the diagonal. See Figure \ref{fig:Ziegelmeierpedexample} for an example of a few snapshots of a nested sequence of Vietoris-Rips simplicial complexes for a point cloud with 18 points as well as the associated barcode and persistence diagram for these data.

\paragraph{Mapper Graph} An alternative applied topology technique is the \emph{mapper graph}, a useful tool for exploration of the structure of a (potentially) high dimensional data set \cite{ZiegelmeierSingh}. To construct a mapper graph from a set of points, we first choose a \emph{filter function} (e.g., a coordinate, density measure, eccentricity, etc.), which assigns a real number to each data point. Then, we select a \emph{cover}, an overlapping collection of sets, over the interval of real values assigned by the filter function. Next, we consider the original set of points with function values in a single set of the cover (the \emph{preimage} of the filter function) and use a clustering algorithm to break these points into distinct groups. Each group becomes a node in the mapper graph, and edges correspond to overlapping groups in the cover. The result encodes a high-dimensional point cloud in a 1-dimensional graph structure such that each node represents a subset of the data points, grouping like points and separating dissimilar points. Obtaining a mapper graph needs a choice of each of these parameters and requires finesse to do so.

\end{document}